\documentclass[12pt]{article}
\usepackage{amssymb}
\usepackage{amsmath}
\usepackage{latexsym}
\usepackage{mathrsfs}

\numberwithin{equation}{section}
\newtheorem{thm}[equation]{Theorem}

\newtheorem{lem}[equation]{Lemma}
\newtheorem{corol}[equation]{Corollary}

\title{Eigenvalues of poly-harmonic operators on variable domains }

\author{Davide Buoso and Pier Domenico Lamberti\footnote{corresponding author: lamberti@math.unipd.it}}

\date{\ }

\begin{document}

\newcommand{\rea}{\mathbb{R}}

\maketitle

%
%
%

\noindent
{\bf Abstract:} We consider a class of eigenvalue problems for poly-harmonic operators,  including Dirichlet and buckling-type
eigenvalue problems.  We prove an analyticity result for the dependence of the symmetric functions of the eigenvalues upon domain perturbations and  compute Hadamard-type formulas for the Frech\'{e}t differentials. We also consider isovolumetric domain perturbations and  characterize the corresponding critical domains for the symmetric functions of the eigenvalues. Finally, we prove that balls are critical domains.
\\

\vspace{11pt}


\vspace{6pt}

\section{Introduction}

Let $n,m\in {\mathbb{N}}_0$ with $0\le m<n$ and  $\Omega $ be a bounded open set  in ${\mathbb{R}}^N$ with smooth boundary. We consider the following eigenvalue problem

\begin{equation}
\label{classic}
{\mathcal{P}}_{nm}:\ \left\{
\begin{array}{ll}
(-\Delta)^n u=\lambda (-\Delta )^mu, & {\rm in}\ \Omega ,\vspace{2mm}\\
u=\frac{\partial u}{\partial \nu}=\dots =\frac{\partial^{n-1} u}{\partial \nu^{n-1}}=0, & {\rm on }\ \partial \Omega ,
\end{array}\right.
\end{equation}
where $\nu$ denotes the unit outer normal to $\partial \Omega$. The case $m=0$ corresponds to the well-known eigenvalue problem for the poly-harmonic operator $(-\Delta )^n$ subject to Dirichlet boundary conditions, while the case $m>0$ represents a buckling-type problem.  These cases include important problems in linear elasticity. For instance, for $N=2$, ${\mathcal{P}}_{10}$ arises in the study of a vibrating membrane stretched in a fixed frame, ${\mathcal{P}}_{20}$
corresponds to the case of a vibrating clamped plate and ${\mathcal{P}}_{21}$ is related to plate buckling. We are mainly interested in the Dirichlet problem ${\mathcal{P}}_{n0}$
and the buckling problem ${\mathcal{P}}_{21}$. However, we prefer to present a unified approach involving all cases.
Problem ${\mathcal{P}}_{nm}$ admits a divergent sequence of positive eigenvalues of finite multiplicity represented as follows
$$
0<\lambda_1[\Omega]\le \lambda_2[\Omega]\le \cdots \le \lambda_j[\Omega ]\le \cdots .
$$
As usual, we agree to repeat each eigenvalue as many times as its multiplicity.

In this paper we are interested in the dependence
of $\lambda_j[\Omega ]$ on $\Omega $. There is a vast literature devoted to domain perturbation problems for elliptic operators. In particular, the cases $n=1,2$ which correspond to the Laplace operator and the bi-harmonic operator respectively, have been intensively studied by many authors.
We refer to Bucur and Buttazzo~\cite{bucbut}, Daners~\cite{dan}, Hale~\cite{hale}, Henry~\cite{henry}, Henrot~\cite{henrot}, Kesavan~\cite{kes} for updated information on this topic.
The case $n>2$ has been much less investigated. However, a renewed general interest in higher order operators has been growing in the last decades as it appears in the extensive monograph by Gazzola, Grunau and Sweers~\cite{gaz} devoted to recent developments in the theory of poly-harmonic operators. As for domain perturbation problems, we refer to the papers \cite{bulahighsh} and \cite{bulahigh} where spectral stability estimates for elliptic operators of arbitrary order are proved. See also the survey paper \cite{bulala}.

Our work is inspired by classical problems in spectral optimization, in particular by the celebrated Rayleigh conjecture (see e.g., \cite{henrot}). Recall that the Rayleigh-Faber-Krahn inequality states that the first eigenvalue of the Laplace operator with Dirichlet boundary conditions (problem ${\mathcal{P}}_{10}$) is minimized by the ball in the class of all bounded domains with a fixed measure. In symbols,
\begin{equation}\label{faber}
\lambda_1[\Omega^* ]\le \lambda_1[\Omega ]\, ,
\end{equation}
where $\Omega ^*$ is a ball with the same measure of $\Omega$. As for the bi-harmonic operator with Dirichlet boundary condition (problem ${\mathcal{P}}_{20}$), inequality (\ref{faber})
was proved by Nadirashvili~\cite{nad} for  $N=2$  and by Ashbaugh and Benguria~\cite{ash} for  $N=3$. See also Mohr~\cite{mohr}.
 Inequality (\ref{faber}) can be proved also for  plate buckling (problem ${\mathcal{P}}_{21}$) under some extra assumptions, see \cite{henrot}. It should be noted that not much is known for higher eigenvalues $\lambda_j[\Omega ]$ for $j>2$. As a corollary of a general result by Buttazzo and Dal Maso~\cite{butdal}, it is known that each eigenvalue of the Dirichlet Laplacian admits a minimizer in the class of all quasi-open sets with fixed measure, contained in a prescribed bounded region.
However, that result says little about the shape of such minimizer.
It is proved in Wolf and Keller~\cite{wolf} that the minimizers of higher eigenvalues in general are not balls and not even unions of balls. Moreover, by looking at the interesting numerical results presented in Oudet~\cite{ou}, one may get the idea that balls are not much relevant in the analysis of higher eigenvalues.

Our main aim is to point out that, despite the above mentioned negative results, balls play a relevant role in the study of isovolumetric perturbations of $\Omega $ for all eigenvalues  $\lambda_j[\Omega ]$ of all problems ${\mathcal{P}}_{nm}$. To do so we shall consider problems ${\mathcal{P}}_{nm}$ on families of
 open sets $\phi (\Omega )$ described by suitable diffeomorphisms $\phi$ defined on a fixed open set  $\Omega$ and we shall study the dependence of $\lambda_j[\phi (\Omega )]$  on $\phi$. As is known, this allows to talk about differentiability and apply calculus in order to find critical eigenvalues with respect to perturbations of $\phi $. One of the main difficulties in the analysis of higher eigenvalues is related to the variation of their multiplicity. This leads to bifurcation
 phenomena which complicate things. For instance, if a fixed open set  $\tilde \phi (\Omega )$ is subject to a perturbation $\phi $ of $\tilde \phi $ then
 a multiple eigenvalue $\lambda_j[\tilde \phi (\Omega )]$ of multiplicity $r$  may split into $r$ simple eigenvalues $\lambda_{j}[ \phi(\Omega )], \dots , \lambda_{j+r-1}[\phi ( \Omega )]$ in such a way that  $\lambda_j[\phi (\Omega )], \dots , \lambda_{j+r-1}[\phi (\Omega )] $ are not differentiable in the variable $\phi$. As for the Laplace operator with Dirichlet or Neumann boundary conditions, it was pointed out  in \cite{lalcri} that in  the case of multiple eigenvalues  it is natural to consider the elementary symmetric functions of the eigenvalues $\lambda_{j}[\phi( \Omega )], \dots , \lambda_{j+r-1}[ \phi(\Omega) ]$.  In this paper, we generalize the results of \cite{lalcri, lala} to all problems ${\mathcal{P}}_{nm}$. Namely, we prove that the elementary symmetric functions of the eigenvalues  $\lambda_{j}[\phi( \Omega )], \dots , \lambda_{j+r-1}[\phi (\Omega )]$ of  ${\mathcal{P}}_{nm}$
 depend real analytically on  $\phi $ (Theorem~\ref{mainthm}) and we prove that if $\tilde\phi (\Omega )$ is a ball then $\tilde \phi$ is  a critical point for such functions under volume constraint (Theorem~\ref{ballthm}). In fact, all critical points $\tilde\phi $  for the symmetric functions of the eigenvalues splitting from an eigenvalue $\lambda $ of multiplicity $r$, can be characterized as those open sets for which the following overdetermined system has a nontrivial solution $(u_1,\dots ,u_r)$
\begin{equation}
\label{classicov}
 \left\{
\begin{array}{lll}
(-\Delta)^n u_i=\lambda (-\Delta )^mu_i, & {\rm in}\ \tilde \phi (\Omega ), & \forall \ i=1, \dots , r ,\vspace{2mm}\\
u_i=\frac{\partial u_i}{\partial \nu}=\dots =\frac{\partial^{n-1} u_i}{\partial \nu^{n-1}}=0, & {\rm on }\ \partial \tilde \phi (\Omega ), & \forall \ i=1, \dots , r,\vspace{2mm}\\
(\frac{\partial^n u_1}{\partial \nu^n})^2+\cdots+(\frac{\partial^n u_r}{\partial \nu^n})^2={\rm const}, & {\rm on }\ \partial \tilde \phi (\Omega ). &
\end{array}\right.
\end{equation}
 Since (\ref{classicov}) is satisfied if $\tilde \phi (\Omega )$ is a ball,  it would be interesting to clarify whether the existence of nontrivial solutions to  (\ref{classicov}) on bounded connected open set $\tilde \phi (\Omega)$ implies that $\tilde \phi (\Omega )$ is a ball.   For $r=1$, $n=1, m=0$ or  $n=2, m=0$  this can be proved by using the celebrated moving plane method under the assumption that the solution $u_1$ does not change sign (see e.g.,  Henry~\cite{he82} for the Laplace operator and Dalmasso~\cite{dalmasso} for the biharmonic operator); for $r=1$, $n=2, m=1$ a different method by Weinberger and Willms leads to the same conclusion (see e.g., \cite{henrot}).

\section{Notation and preliminaries}

Let $n\in {\mathbb{N}}$ and  $\Omega$ be a bounded open set in ${\mathbb{R}^N}$.
By $W^{n,2}(\Omega )$ we denote the Sobolev space of all functions in $L^2(\Omega )$ which admit weak derivatives in $L^2(\Omega )$ up to order $n$. By $W^{n,2}_0(\Omega )$ we denote the closure in $W^{n,2}(\Omega )$ of the space of $C^{\infty }$-functions with compact support in $\Omega $.
We consider the weak formulation of problem (\ref{classic}). To do so, for any $m\in {\mathbb{N}}_0$ with  $0\le m \le n$,  we consider the poly-harmonic operator $\Delta ^m$ as the operator from $W^{n,2}_0(\Omega )$ to its dual $(W^{n,2}_0(\Omega ))'$ which takes any $u\in W^{n,2}_0(\Omega )$ to the functional
$\Delta ^m[u]$ defined by
\begin{equation}\label{pol1cla}
\Delta^{2s}[u][\varphi ]=\int_{\Omega }\Delta^su\Delta^s\varphi dx,\ \ \ \forall\varphi  \in W^{n,2}_0(\Omega ),
\end{equation}
if $m=2s$ and
\begin{equation}\label{pol2cla}
\Delta^{2s+1}[u][\varphi ]=-\int_{\Omega }\nabla (\Delta^su)\cdot \nabla (\Delta^s\varphi) dx,\ \ \ \forall\varphi  \in W^{n,2}_0(\Omega ),
\end{equation}
if $m=2s+1$, where $s\in {\mathbb{N}}_0$. Thus,  the weak formulation of the classic problem (\ref{classic}) reads
\begin{equation}\label{weak}
(-\Delta )^n[u][\varphi ]=\lambda (-\Delta )^m[u][\varphi ],\ \ \forall \varphi \in W^{n,2}_0(\Omega ).
\end{equation}

By the Poincar\'{e} inequality, it follows that the quadratic form defined by $(-\Delta )^n[u][u] $ for all $u\in W^{n,2}_0(\Omega )$ is coercive in
$W^{n,2}_0(\Omega )$, hence the operator $(-\Delta )^n $ is a linear homeomorphism from $W^{n,2}_0(\Omega )$ onto  $(W^{n,2}_0(\Omega ))'$.
Thus the equation (\ref{weak}) is equivalent to the equation
$
(-\Delta )^{-n}\circ (-\Delta )^m [u] = \lambda^{-1} u,
$
where $(-\Delta )^{-n}$ denotes the inverse of $(-\Delta )^{n}$. It is convenient  to endow the space $W^{n,2}_0(\Omega )$ with the scalar product defined by
\begin{equation}\label{scalar}
<u_1,u_2>_{n}=(-\Delta )^n[u_1][u_2],
\end{equation} for all $u_1,u_2\in W^{n,2}_0(\Omega )$. The norm induced by this scalar product is equivalent to the standard Sobolev norm. In the sequel, unless otherwise indicated, we shall think of  $W^{n,2}_0(\Omega )$ as a Hilbert space equipped with the scalar product (\ref{scalar}). This  allows to give a straightforward proof of the following

\begin{lem}\label{comp} Let $\Omega $ be a bounded open set in ${\mathbb{R}}^N$ and $m,n\in {\mathbb{N}} $ with $0\le m<n$.
 The operator  $S_{\Omega }\equiv (-\Delta )^{-n}\circ (-\Delta )^m $ is a non-negative self-adjoint compact  operator in the Hilbert space  $W^{n,2}_0(\Omega )$. The spectrum of $S_{\Omega }$ is discrete and consists of a decreasing sequence of positive eigenvalues of finite multiplicity converging to zero. Moreover, the equation $S_{\Omega } u=\mu u$ is satisfied for some $u\in W^{n,2}_0(\Omega )$, $\mu >0$ if and only if
 equation (\ref{weak}) is satisfied with $\lambda =\mu^{-1}$.
\end{lem}

{\bf Proof.} The equality  $
<S_{\Omega } u_1,u_2>_n=(-\Delta)^m[u_1][u_2]$, for all $u_1,u_2\in W^{n,2}_0(\Omega )$ and the symmetry of the operator $(-\Delta )^m$ implies that $ S_{\Omega }$ is a self-adjoint operator. Since $\Omega$ is bounded and $m<n$, the space $W^{n,2}_0(\Omega )$ is compactly embedded into $W^{m,2}_0(\Omega )$. This implies that the operator $(-\Delta )^m$ is a compact operator from the space $W^{n,2}_0(\Omega )$ to its dual. The rest of the proof is trivial. \hfill $\Box$\\

By Lemma \ref{comp} and standard spectral theory we deduce the following

\begin{corol} Let $\Omega $ be a bounded open set in ${\mathbb{R}}^N$ and $m,n\in {\mathbb{N}} $ with $0\le m<n$. The eigenvalues of problem (\ref{weak}) are positive, have finite multiplicity and can be represented as a  non-decreasing divergent sequence $\lambda_j[\Omega ]$, $j\in {\mathbb{N}}$ where each eigenvalue is repeated according to its multiplicity. Moreover,
\begin{equation}
\label{minmax}
\lambda_j[\Omega ]\equiv \lambda_j^{n,m}[\Omega ]=\min_{\substack{E\subset W^{n,2}_0(\Omega ) \\     {\rm dim } E=j}} \max_{\substack{ u\in E \\  u\ne 0}} R_{nm}[u],
\end{equation}
for all $j\in {\mathbb{N}}$, where $R_{nm}[u]$ is the Rayleigh quotient defined by
$$
R_{nm}[u]=\left\{\begin{array}{lll} \frac{  \int_{\Omega }| \Delta^ru|^2 dx  }{ \int_{\Omega }|\Delta^su|^2 dx    },& \ \ {\rm if}\ n=2r,&\ m=2s,\vspace{2.5mm} \\
\frac{  \int_{\Omega }|\Delta^ru|^2 dx  }{ \int_{\Omega }|\nabla \Delta^su|^2 dx    },& \ \ {\rm if}\ n=2r,&\ m=2s+1,\vspace{2.5mm} \\
\frac{  \int_{\Omega }|\nabla \Delta^ru|^2 dx  }{ \int_{\Omega }|\Delta^su|^2 dx    },& \ \ {\rm if}\ n=2r+1,&\ m=2s,\vspace{2.5mm} \\
\frac{  \int_{\Omega }|\nabla \Delta^ru|^2 dx  }{ \int_{\Omega }|\nabla\Delta^su|^2 dx    },& \ \ {\rm if}\ n=2r+1,&\ m=2s+1.    \end{array}  \right.
$$
\end{corol}

Clearly, the eigenvalues $\lambda _j^{n,m}[\Omega ]$  depend on $n,m$. However,
for the sake of simplicity, we shall write $\lambda _j[\Omega ]$ instead of $\lambda _j^{n,m}[\Omega ]$.

\section{Analyticity results}

Let $\Omega $ be a bounded open set in ${\mathbb{R}}^N$ of class $C^1$.  In the sequel, we shall consider problem (\ref{weak}) in a family of open sets  parameterized by suitable diffeomorphisms $\phi $
defined on $\Omega $. Namely, we set
$$
{\mathcal{A}}_{\Omega }^n=\biggl\{\phi\in C^n_b(\Omega\, ; {\mathbb{R}}^N ):\ \inf_{\substack{x_1,x_2\in \Omega \\ x_1\ne x_2}}\frac{|\phi(x_1)-\phi(x_2)|}{|x_1-x_2|}>0 \biggr\},
$$
where $C^n_b(\Omega\, ; {\mathbb{R}}^N )$ denotes the space of all functions from $\Omega $ to ${\mathbb{R}}^N$ of class $C^n$, with bounded derivatives up to order
  $n$.  Note that if $\phi \in {\mathcal{A}}_{\Omega }^n$ then $\phi $ is injective, Lipschitz continuous and $\inf_{\Omega }|{\rm det }\nabla \phi |>0$. Moreover, $\phi (\Omega )$ is a bounded open set and the inverse map $\phi^{(-1)}$ belongs to  ${\mathcal{A}}_{\phi(\Omega )}^n$.  Thus it is natural to consider problem (\ref{weak}) on $\phi (\Omega )$ and study  the dependence of $\lambda_j[\phi (\Omega )]$ on $\phi \in {\mathcal{A}}_{\Omega }^n$. To do so, we endow the space $C^n_b(\Omega\, ; {\mathbb{R}}^N )$ with its usual norm $\| f \|_{C^n_b(\Omega\, ;{\mathbb{R}}^N )}=\sup_{|\alpha|\le n,\ x\in\Omega  } |D^{\alpha }f(x)|$. Note that ${\mathcal{A}}_{\Omega }^n$ is an open set in  $C^n_b(\Omega\, ;{\mathbb{R}}^N )$, see \cite[Lemma~3.11]{lala}.
    Thus, it makes sense to  study differentiability and analyticity properties of the maps $\phi \mapsto \lambda_j[\phi (\Omega )]$ defined for $\phi \in {\mathcal{A}}_{\Omega }^n$.
   For simplicity, we write  $\lambda_j[\phi ]$ instead of $\lambda_j[\phi (\Omega )]$.
   As in \cite{lala}, we fix a finite set of indexes $F\subset \mathbb{N}$
and we consider those maps $\phi\in {\mathcal{A}}^n_{\Omega }$ for which the eigenvalues
with indexes in $F$  do not coincide with eigenvalues with indexes not
in $F$; namely we set
$$
{\mathcal { A}}^{n,m}_{F, \Omega }= \left\{\phi \in {\mathcal { A}}^n_{\Omega }:\
\lambda_j[\phi ]\ne \lambda_l[\phi],\ \forall\  j\in F,\,   l\in \mathbb{N}\setminus F
\right\}.
$$
It is also convenient to consider those maps $\phi \in {\mathcal { A}}^{n,m}_{F, \Omega } $ such that all the eigenvalues with index in $F$
 coincide and set
$$
\Theta_{F, \Omega }^{n,m} = \left\{\phi\in {\mathcal { A}}_{F, \Omega }^{n,m}:\ \lambda_{j_1}[\phi ]
=\lambda_{j_2}[\phi ],\, \
\forall\ j_1,j_2\in F  \right\} .
$$

 For $\phi \in {\mathcal { A}}^{n,m}_{F, \Omega }$, the elementary symmetric functions of the eigenvalues with index in $F$ are defined by
\begin{equation}
\label{sym1}
\Lambda_{F,h}[\phi ]=\sum_{ \substack{ j_1,\dots ,j_h\in F\\ j_1<\dots <j_h} }
\lambda_{j_1}[\phi ]\cdots \lambda_{j_h}[\phi ],\ \ \ h=1,\dots , |F|.
\end{equation}

The main  result of this section is the following generalization to poly-harmonic operators on smooth domains of the results in \cite[\S 3]{lala} concerning the Dirichlet Laplacian on rough domains.

\begin{thm}\label{mainthm}Let $\Omega $ be a bounded open set in ${\mathbb{R}}^N$ of class $C^1$, $n,m\in {\mathbb{N}}_0$ with $0\le m <n$,  and  $F$ be a finite set in  ${\mathbb{N}}$. The set ${\mathcal { A}}^{n,m}_{F, \Omega }$ is open in
$C^n_b(\Omega\, ; {\mathbb{R}}^N )$ and the real-valued maps which takes $\phi\in {\mathcal { A}}^{n,m}_{F, \Omega } $ to $ \Lambda_{F,h}[\phi ]$ are real-analytic on  ${\mathcal { A}}^{n,m}_{F, \Omega }$ for all $h=1,\dots , |F|$. Moreover, if $\tilde \phi\in \Theta^{n,m}_{F, \Omega }  $ is such that the eigenvalues $\lambda_j[\tilde \phi]$ assume the common value $\lambda_F[\tilde \phi ]$ for all $j\in F$, and  $\tilde \phi (\Omega )$ is of class $C^{2n}$ then  the Frech\'{e}t differential of the map $\Lambda_{F,h}$ at the point $\tilde\phi $ is delivered by the formula
\begin{equation}
		\label{derivn}
			d|_{\phi=\tilde{\phi}}\Lambda_{F,h}[\psi]=-\lambda_F^h[\tilde{\phi}]\binom{|F|-1}{h-1}
			\sum_{l=1}^{|F|}{\int_{\partial \tilde{\phi}(\Omega)}}  \left(  \frac{\partial^n v_l}{\partial \nu^n}  \right)^2(\psi\circ\tilde{\phi}^{(-1)})\cdot\nu d\sigma,
		\end{equation}
		for all $\psi\in C^{n}_b(\Omega ; {\mathbb{R}}^N)$, where $\{v_l\}_{l\in F}$ is an orthonormal basis in $W^{n,2}_0(\tilde \phi (\Omega ))$ (with respect to the scalar product (\ref{scalar})) of the eigenspace associated with $\lambda_F[\tilde \phi]$, and $\nu $ denotes the unit outer normal to $\partial\tilde \phi (\Omega )$.
\end{thm}

Note that formula (\ref{derivn}) is a generalization of the celebrated  Hadamard formula, see Grinfeld~\cite{grin} for a recent paper on this topic;
see also Ortega and Zuazua~\cite{zua} for the analysis of associated bifurcation phenomena concerning multiple eigenvalues of the biharmonic operator subject to Dirichlet boundary conditions.

In order to prove Theorem \ref{mainthm} we consider  equation (\ref{weak}) on $\phi (\Omega )$ and pull it back to $\Omega $.
 Namely, we consider the equation
\begin{equation}\label{weakvv}
(-\Delta)^n[v][\psi ]=\lambda (-\Delta)^m[v][\psi ],\ \ \forall \ \psi \in W^{n,2}_0(\phi(\Omega )),
\end{equation}
in the unknowns $v\in W^{n,2}_0(\phi (\Omega ))$, $\lambda\in ]0,\infty [$. Recall that the pull-back to $\Omega $ of the classic Laplace operator on $\phi (\Omega)$ is defined by
\begin{equation}\label{lapphi}
\Delta_{\phi }u=(\Delta (u\circ \phi^{(-1)}) )\circ \phi
\end{equation}
for all $u\in W^{2,1}_{loc}(\Omega )$, $\phi \in {\mathcal{A}}_{\Omega }^2$. The operator $\Delta_{\phi}$ is in fact the Laplace-Beltrami operator associated with the change of variables defined by $\phi$.  Note that
\begin{equation}\label{polybelt}
\Delta_{\phi }^su =(\Delta^s (u\circ \phi^{(-1)}) )\circ \phi
\end{equation}
for all $u\in W^{2s,1}_{loc}(\Omega)$, $\phi \in {\mathcal{A}}_{\Omega }^{2s}$.
For any $0\le m \le n$,
 the operator $\Delta_{\phi }^m$ can be considered as the operator acting from
$ W^{n,2}_0(\Omega ) $ to its dual, which takes any $u\in W^{n,2}_0(\Omega )$ to the functional $\Delta_{\phi}^n[u]$ defined by
$$
\Delta_{\phi}^m[u][\varphi]=\Delta^m [u\circ \phi^{(-1)}][\varphi\circ \phi^{(-1)}],
$$
for all $\varphi\in W^{n,2}_0(\Omega ) $. More precisely,  if $m=2s$,  $s\in {\mathbb{N}}_0$ then
\begin{equation}\label{pol1}
\Delta^{2s}_{\phi }[u][\varphi ]=     \int_{\Omega }\Delta^s_{\phi} u\Delta^s_{\phi }\varphi |{\rm det}\nabla \phi|dx,
\end{equation}
for all $\varphi  \in W^{n,2}_0(\Omega )$.  If $m=2s+1$,  $s\in {\mathbb{N}}_0$ then
\begin{equation}\label{pol2}
-\Delta^{2s+1}_{\phi }[u][\varphi ]=\int_{\Omega }\nabla (\Delta^s_{\phi}u)\ (\nabla\phi)^{-1}(\nabla\phi)^{-t} \nabla^t (\Delta^s_{\phi }\varphi)|{\rm det}\nabla \phi| dx,
\end{equation}
for all  $\varphi \in W^{n,2}_0(\Omega )$, where $(\nabla\phi )^{-1}$ denotes the inverse of the Jacobian matrix of $\phi $ and $(\nabla\phi )^{-t}$ the transpose of $(\nabla\phi )^{-1}$.
 Note that the map from $W^{n,2}_0(\Omega )$ to  $W^{n,2}_0(\phi (\Omega ))$ which maps $u$ to $u\circ \phi^{(-1)} $ for all $u\in  W^{n,2}_0(\Omega )$ is a linear homeomorphism. Hence, equation (\ref{weakvv}) is equivalent to
\begin{equation}
(-\Delta_{\phi })^n[u][\varphi ]=\lambda (-\Delta_{\phi })^m[u][\varphi ],\ \ \forall \ \varphi \in W^{n,2}_0(\Omega )
\end{equation}
where $u=v\circ\phi$. It is also natural to pull-back the scalar product of $W^{n,2}_0(\phi (\Omega ))$ to $\Omega $ by setting
\begin{equation}
\label{sca}
<u_1,u_2>_{n,\phi}=<u_1\circ \phi^{(-1)}, u_2\circ \phi^{(-1)} >_n
\end{equation}
for all $u_1,u_2 \in W^{n,2}_0(\Omega )$, where $<\cdot , \cdot >_n$ is the scalar product in $W^{n,2}_0(\phi (\Omega ))$ defined by (\ref{scalar}). By $W^{n,2}_{0,\phi }(\Omega )$ we denote the Hilbert space $W^{n,2}_0(\Omega )$ endowed with the scalar product $<\cdot ,\cdot >_{n, \phi }$. It turns out that the operator $S_{\phi (\Omega )}$ defined in Lemma \ref{comp} is unitarily equivalent to the operator $T_{\phi }$ defined on $W^{n,2}_{0, \phi}(\Omega )$  by
\begin{equation}\label{tphi}
T_{\phi }= (-\Delta_{\phi })^{-n}\circ (-\Delta_{\phi })^m.
\end{equation}

Thus we can prove the following lemma where ${\mathcal{L}}(W^{n,2}_{0}(\Omega ))$  denotes the space of linear bounded operators from $W^{n,2}_{0}(\Omega )$ to itself and and ${\mathcal{B}}_s(W^{n,2}_{0}(\Omega ))$ denotes the space of bilinear forms on $W^{n,2}_{0}(\Omega )$ (both spaces are equipped with their usual norms).

\begin{lem}\label{changelemma} Let $\Omega $ be a bounded open set in ${\mathbb{R}}^N$ of class $C^1$, $n,m\in {\mathbb{N}}_0$, $0\le m <n$.  The operator $T_{\phi}$ defined in (\ref{tphi}) is  non-negative self-adjoint and compact  on the Hilbert space $W^{n,2}_{0, \phi}(\Omega )$. The equation (\ref{weakvv}) is satisfied for some $v\in W^{n,2}_0(\phi(\Omega ))$ if and only if  the equation $T_{\phi } u=\mu u$ is satisfied with
$u=v\circ \phi $ and $\mu =\lambda^{-1}$. Moreover, the map from ${\mathcal{A}}_{\Omega }^n$ to ${\mathcal{L}}(W^{n,2}_{0}(\Omega ))\times {\mathcal{B}}_s(W^{n,2}_{0}(\Omega ))$  which takes  $\phi \in {\mathcal{A}}_{\Omega }^n$ to $(T_{\phi }, <\cdot ,\cdot>_{n,\phi})$ is real-analytic.
\end{lem}

{\bf Proof.} Since  the operator $T_{\phi }$ is unitarily equivalent to the operator  $S_{\phi (\Omega )}$,
 the first part of the lemma immediately follows by Lemma~\ref{comp}. In order to prove the real-analytic dependence of $T_{\phi}$ upon $\phi$, we note that by standard calculus
\begin{equation}\label{change}
 \Delta _{\phi }u=\sum_{r,s,i=1}^N\left(\frac{\partial^2 u}{\partial x_r\partial x_s}\sigma_{ri}\sigma_{si}+ \frac{\partial u}{\partial x_r}\frac{\partial \sigma _{ri}}{\partial x_s}\sigma_{si}\right)
\end{equation}
 for all $u\in W^{2,2}(\Omega )$, where $\sigma =(\nabla \phi )^{-1}$ (see also \cite[Proposition~3.1]{lala}). By formula (\ref{change}),  it follows that the map from ${\mathcal{A}}^n_{\Omega }\times W^{n,2}(\Omega )$ to $W^{n-2,2}(\Omega )$ which takes $(\phi ,u) \in {\mathcal{A}}^n_{\Omega }$ to $\Delta_{\phi }u$ is real-analytic. Thus also the maps from ${\mathcal{A}}^n_{\Omega } \times W^{n,2}(\Omega ) $ to $L^2(\Omega )$ which take $(\phi ,u) \in {\mathcal{A}}^n_{\Omega }$ to $\Delta_{\phi }^{s}u$ for all $s\in {\mathbb{N}}_0$ with $0\le s\le n/2$,  are real-analytic since they are compositions of real-analytic maps.
This, combined with formulas (\ref{pol1}) and (\ref{pol2}), implies the real-analytic dependence of $T_{\phi}$ and $<\cdot ,\cdot >_{n,\phi}$ upon $\phi$.\hfill $\Box$ \\

{\bf Proof of Theorem \ref{mainthm}.} We denote by $\mu_j[\phi]$, $j \in {\mathbb{N}}$, the eigenvalues of $T_{\phi}$.
By Lemma~\ref{changelemma}, $\mu_j[\phi]=\lambda_j^{-1}[\phi ]$ for all $j\in {\mathbb{N}}$, hence the set  ${\mathcal { A}}^{n,m}_{F, \Omega }$  coincides
with the set
$\{\phi\in {\mathcal{A}}_{\Omega }^n:\
\mu_j[\phi]\ne \mu_l[\phi],\  \forall j\in F,\   l\in \mathbb{N}\setminus F
\} $. By Lemma~\ref{changelemma},
 $T_{\phi }$
is self-adjoint with respect to the scalar product
$<\cdot ,\cdot >_{n,\phi} $ and both $T_{\phi }$ and $<\cdot ,\cdot >_{n,\phi} $ depend real-analytically on $\phi$.
Thus, by applying \cite[Thm.~2.30]{lala},
it follows that
${\mathcal { A}}^{n,m}_{F, \Omega }$ is an open set in $C^n_b(\Omega\, ; {\mathbb{R}}^N )$  and the
functions which take $\phi\in {\mathcal { A}}^{n,m}_{F, \Omega } $  to
\begin{equation}
\label{sym3}
\Gamma_{F,h}[\phi ]=\sum_{\substack{j_1,\dots ,j_h\in F\\ j_1<\dots <j_h}}
\mu_{j_1}[\phi ]\cdots \mu_{j_h}[\phi ]
\end{equation}
are real-analytic for all $h=1, \dots , |F|$.
Since
\begin{equation}
\label{sym4}
\Lambda_{F,h}[\phi ] = \frac{\Gamma_{F,|F|-h}[\phi ] }{\Gamma_{F,|F|}[\phi ]},
\end{equation}
for all $h=1, \dots ,|F|$, where $\Gamma_{F,0}[\phi ]= 1  $,
it follows that $\Lambda_{F,h}[\phi ]$ depends real-analytically
on $\phi\in {\mathcal { A}}^{n,m}_{F, \Omega }$.

It remains to  prove formula (\ref{derivn}).  Let $\tilde \phi \in \Theta_{F, \Omega }^{n,m} $, $\lambda_F[\tilde \phi ]$ and $\{  v_l\}_{l\in F}$ be as in the statement. We set $ u_l= v_l\circ \tilde \phi $ for all $l\in F$ and we note that $\{  u_l\}_{l\in F}  $ is an orthonormal basis in $W^{n,2}_{0,\tilde \phi }(\Omega )$ for the eigenspace corresponding to the eigenvalue $\lambda ^{-1}_F[\tilde \phi ]$ of the operator $T_{\tilde \phi }$.
By
\cite[Thm.~2.30]{lala}, it follows that
\begin{equation}
\label{sym5}
{\rm d|}_{\phi =\tilde\phi}\Gamma_{F,h}[\psi]=
\lambda_{F}^{1-h}[\tilde \phi]    {|F|-1\choose h-1}
\sum_{l\in F}
<{\rm d|}_{\phi =\tilde \phi}T_{\phi }[\psi][u_l], u_l>_{n,\tilde \phi}
\end{equation}
for all  $\psi\in C^n_b( \Omega\, ;{\mathbb{R}}^N) $.   Note that by standard regularity theory (see e.g., Agmon~\cite[Thm.~9.8]{ag}) $v_l\in W^{2n,2}(\tilde \phi (\Omega ))$ for all $l\in F$.

By standard calculus, equalities (\ref{scalar}), (\ref{sca}), Theorem \ref{beldif}, by observing that
$\frac{\partial ^m v_l}{\partial \nu ^m}=0$  on $\partial \tilde\phi (\Omega )$ and $(-\Delta )^n v_l=\lambda_F[\tilde \phi](-\Delta )^m v_l$,
we have

\begin{eqnarray}\label{sym6}\lefteqn{
<{\rm d|}_{\phi =\tilde\phi } T_{\phi }[\psi][u_l], u_l>_{n,\tilde \phi}}\nonumber \\ & &
= ({\rm d|}_{\phi =\tilde\phi } (-\Delta _{\phi})^m[\psi])[u_l][u_l] -\lambda_F^{-1}[\tilde\phi] ({\rm d|}_{\phi =\tilde\phi } (-\Delta _{\phi})^n[\psi])[u_l][u_l]\nonumber  \\ & &
= -\int_{\partial \tilde \phi (\Omega )} \left(\frac{\partial ^mv_l}{\partial \nu ^m}\right)^2\zeta \cdot \nu  d\sigma
-2\int_{\partial \tilde \phi (\Omega )}(-\Delta )^mv_l\nabla v_l\cdot \zeta d\sigma \nonumber \\ & &
+ \lambda_F^{-1}[\tilde \phi ]\int_{\partial \tilde \phi (\Omega )} \left(\frac{\partial ^nv_l}{\partial \nu ^n}\right)^2 \zeta \cdot \nu d\sigma
+2\lambda_F^{-1}[\tilde \phi ]\int_{\partial \tilde \phi (\Omega )}(-\Delta )^nv_l\nabla v_l\cdot \zeta d\sigma \nonumber  \\
& & = \lambda_F^{-1}[\tilde \phi ]\int_{\partial \tilde \phi (\Omega )} \left(\frac{\partial ^nv_l}{\partial \nu ^n}\right)^2 \zeta \cdot \nu d\sigma ,
\end{eqnarray}
where we have set $\zeta =\psi \circ \tilde \phi^{(-1)}$.
Formula (\ref{derivn}) easily follows by combining formulas (\ref{sym4})-(\ref{sym6}). \hfill $\Box$\\

\section{Isovolumetric perturbations}

Consider the following extremum problems for the symmetric functions of the eigenvalues
\begin{equation}\label{min}
\min_{V[\phi ]={\rm const}}\Lambda_{F,h}[\phi ]\ \ \ {\rm or}\ \ \ \max_{   V[\phi ]={\rm const}}\Lambda_{F,h}[\phi],
\end{equation}
where $V[\phi ]$ denotes the $N$-dimensional Lebesgue measure of $\phi (\Omega )$.

Note that if $\tilde \phi \in {\mathcal{A}}_{\Omega }^n$ is a minimizer or maximizer in (\ref{min}) then $\tilde \phi $ is a critical domain transformation for the map $\phi\mapsto  \Lambda_{F,h}[\phi ]$ subject to volume constraint, {\rm i.e.,}
\begin{equation}\label{inc}
{\rm Ker\ d}|_{\phi =\tilde \phi }V \subset {\rm Ker\ d}|_{\phi =\tilde \phi } \Lambda_{F,h},
\end{equation}
where $V$ is the real valued function defined on ${\mathcal{A}}_{\Omega }^n$ which takes $\phi \in {\mathcal{A}}_{\Omega }^n$ to $V[\phi ]$.

The following theorem provides a characterization of all critical domain transformations $\phi$. See \cite{lalcri}
for the case of the Dirichlet and Neumann Laplacians.

\begin{thm}\label{car}Let $\Omega $ be a bounded open set in ${\mathbb{R}}^N$ of class $C^1$, $n,m\in {\mathbb{N}}_0$ with $0\le m<n$,  and $F$ be a finite subset of ${\mathbb{N}}$.
Assume that $\tilde \phi\in \Theta_{F, \Omega }^{n,m} $ is such that $\tilde \phi (\Omega )$ is of class $C^{2n}$ and that the eigenvalues $\lambda_j[\tilde \phi ]$ have the common value $\lambda_F[\tilde \phi ]$
for all $j\in F$. Let $\{  v_l\}_{l\in F}$ be an orthornormal basis in $W^{n,2}_0(\tilde \phi (\Omega ))$ of the eigenspace corresponding to $\lambda_F[\tilde \phi ]$. Then  $\tilde \phi$
is a critical domain transformation  for any of the functions $\Lambda_{F,h}$, $h=1,\dots , |F|$,  with  volume constraint
if and only if  there exists $C \in {\mathbb{R}}$ such that
\begin{equation}\label{ov1}
\sum_{l\in F}
 \biggl|\frac{\partial^n  v_l}{\partial \nu^n}\biggr|^2 =C,\ \ {\rm on}\ \partial\tilde \phi (\Omega )\, .
\end{equation}
\end{thm}

{\bf Proof.}  Note that $V[\phi ]=\int_{\Omega }|{\rm det }\nabla \phi |dx$, hence by formula (\ref{der4}) it follows that
\begin{equation}\label{dev}
{\rm d|}_{\phi =\tilde \phi }V[\psi ]= \int_{\partial \tilde \phi (\Omega )}(\psi\circ\tilde \phi ^{(-1)})\cdot \nu d\sigma ,
\end{equation}
for all $\psi \in C^n_b(\Omega\, ; {\mathbb{R}}^N)$.  The proof of (\ref{ov1}) follows  immediately  by formulas (\ref{derivn}) and (\ref{dev}), and  by observing that condition (\ref{inc}) is satisfied if and only if there exists  $c\in {\mathbb{R}}$ (a Lagrange multiplier) such that
\begin{equation}\label{lag}
{\rm d}|_{\phi =\tilde \phi } \Lambda_{F,h} =c\, {\rm d}|_{\phi =\tilde \phi }V.
\end{equation}
\hfill $\Box$

Finally, we can prove the following

\begin{thm}\label{ballthm}
Let the same assumptions of Theorem~\ref{car} hold. If $\tilde \phi (\Omega )$ is a ball then condition (\ref{ov1}) is satisfied.
\end{thm}

{\bf Proof. }  Without any loss of generality, we assume that $\tilde \phi (\Omega )$ is a ball with radius $R$ centered at zero.  By the rotation
invariance of the Laplace operator,  $\{ v_{l}\circ A\}_{l\in F}$   is an orthonormal basis of the eigenspace  corresponding to $\lambda_F[\tilde \phi ]  $ for all $A\in O_{N}({\mathbb{R}})$ where  $O_{N}({\mathbb{R}})$ denotes the group of orthogonal linear transformations in ${\mathbb{R}}^N$.
 Since both $\{ v_l\}_{l\in F}$
and $\{  v_{l}\circ A\}_{l\in F}$ are orthonormal bases of the same space, it follows that
$
\sum_{l=1}^{|F|} { v}_{l}^{2}\circ
A=\sum_{l=1}^{|F|} { v}_{l}^{2}\, ,
$
 for all $A\in O_{n}({\mathbb{R}})$. Thus $\sum_{l=1}^{|F|}{ v}_{l}^{2}$ is a radial function. By differentiating with respect to the radial coordinate $r$, by the Leibniz formula and  by recalling that all derivatives up to order $n-1$ of the eigenfunctions  vanish at the boundary of $\tilde \phi (\Omega )$, we obtain that
 \begin{equation}\label{lei}
 \frac{\partial^{2n} v_l^2 }{\partial r ^{2n} }\biggl|_{r=R}\biggr.=\sum_{k=0}^{2n}\binom{2n}{k}\left(  \frac{\partial^{k}  v_l }{\partial r ^{k} }\frac{\partial^{2n-k}  v_l }{\partial r ^{2n-k} }\right)\biggl|_{r=R}\biggr.=\binom{2n}{n} \left(\frac{\partial^{n}  v_l }{\partial r ^{n} } \right)^2 \biggl|_{r=R}\biggr. \, .
 \end{equation}
Since $\sum_{l\in F} \frac{\partial^{2n}  v_l^2 }{\partial r ^{2n} }$ is a radial function, then by formula (\ref{lei}) we conclude  that the
$\sum_{l\in F}(\frac{\partial^{n}  v_l }{\partial \nu  ^{n} })^2$ is constant on $\partial \tilde \phi (\Omega )$. \hfill $\Box$\\

It would be interesting to clarify whether balls are local minimizers or maximizers for the eigenvalues or their symmetric functions. With regard to this, we mention that  it is proved in Wolf and Keller~\cite[Thm.~8.3]{wolf} that the circular disk in the plane is a local minimizer for the third eigenvalue of the Dirichlet Laplacian.

\section{A formula for the Frech\'{e}t differential of the `poly-Laplace-Beltrami' operator}

In this section we prove Theorem~\ref{beldif}  which has its own interest since it provides an explicit formula for the Frech\'{e}t differential  with respect to $\phi $ of the weak `poly-Laplace-Beltrami' operator  $\Delta_{\phi}^n$ defined by (\ref{pol1}), (\ref{pol2}).
That formula  has been used in the proof of  (\ref{derivn}).

\begin{lem} Let $\Omega $ be a bounded open set in ${\mathbb{R}}^N$ of class $C^1$, $s\in {\mathbb{N}}$,  $u_1\in L^{2}(\Omega)$, $u_2\in W^{2s,2}_0(\Omega)$.
Let $\tilde \phi \in {\mathcal{A}}_{\Omega}^{2s}$ and $v_i=u_i\circ \tilde \phi ^{(-1)}$, $i=1,2$. Assume that $\tilde \phi (\Omega )$ is of class $C^1$ and that $v_1\in W^{2s,2}(\tilde\phi(\Omega))$, $v_2\in W^{2s+1,2}(\tilde\phi(\Omega))$. Then
\begin{eqnarray}\label{applem}
\lefteqn{\int_{\Omega} u_1{\rm d|}_{\phi =\tilde \phi } \Delta^s_{\phi }[\psi]u_2|{\rm det }\nabla\tilde  \phi | dx}\nonumber \\ & & \qquad\
=\int_{\tilde \phi (\Omega )}(v_1\nabla\Delta^sv_2-\Delta ^sv_1\nabla v_2)\cdot \zeta dy
-\int_{\partial \tilde\phi (\Omega )}v_1\Delta^sv_2\zeta \cdot \nu d\sigma ,
\end{eqnarray}
for all $\psi \in C^{2s}_b(\Omega\, ; {\mathbb{R}}^N)$, where $\zeta =\psi \circ \tilde\phi^{(-1)}$.
\end{lem}

{\bf Proof.} First, we recall the following formula from \cite[Lemma~3.42]{lalamult} which holds for any $u\in W^{2,2}(\Omega )$:
\begin{equation}\label{diff1}
({\rm d|}_{\phi =\tilde \phi } \Delta_{\phi }[\psi]u)\circ \tilde \phi^{(-1)}=-2\sum_{i,j=1}^N
\frac{\partial^2 (u\circ \tilde\phi^{(-1)})}{\partial y_i\partial y_j}\frac{\partial \zeta_j}{\partial y_i}
-\sum_{j=1}^N\frac{\partial  (u\circ \tilde\phi^{(-1)})}{\partial y_j}\Delta \zeta _j
\, .
\end{equation}

We observe that
\begin{equation}\label{diff2}d|_{\phi=\tilde{\phi}}\Delta_{\phi}^s[\psi]=\sum_{ \substack{h,k=0\\h+k=s-1}}^{s-1}\Delta_{\tilde{\phi}}^h\circ(d|_{\phi=\tilde{\phi}}\Delta_{\phi}[\psi])
\circ\Delta_{\tilde{\phi}}^k\, ,
\end{equation}
	
By formulas (\ref{diff1}) and (\ref{diff2}), by changing variables in integrals and integrating by parts, we obtain
\begin{eqnarray}\label{ppp0}\lefteqn{
\int_{\Omega} u_1{\rm d|}_{\phi =\tilde \phi } \Delta^s_{\phi }[\psi]u_2|{\rm det }\nabla\tilde \phi |dx } \nonumber \\
\qquad & & =  -\!\!\!\!\!\sum_{ \substack{h,k=0\\h+k=s-1}}^{s-1} \int_{\tilde \phi (\Omega )}\Delta ^hv_1   \left( 2 \sum_{i,j=1}^N \frac{\partial^2  \Delta^kv_2  }{\partial y_i\partial y_j}\frac{\partial \zeta_j}{\partial y_i} +\sum_{j=1}^N\frac{\partial  \Delta^kv_2 }{\partial y_j}\Delta \zeta _j    \right )dy\nonumber \\
\qquad & & =  \sum_{ \substack{h,k=0\\h+k=s-1}}^{s-1} \int_{\tilde \phi (\Omega )} \sum_{i,j=1}^N\left[ \frac{\partial \Delta^hv_1}{\partial y_i }\frac{\partial \Delta^kv_2}{\partial y_j }\left( \frac{\partial \zeta_j}{\partial y_i}   +\frac{\partial \zeta_i}{\partial y_j}   \right)\right] \nonumber  \\
\qquad & &\left.  -(
\Delta^hv_1\Delta^{k+1}v_2
 + \nabla \Delta^h v_1\nabla \Delta^kv_2
){\rm div} \zeta \right. dy,
\end{eqnarray}
see also \cite[Formula (3.45)]{lalamult}.
Moreover, integrating by parts yields

\begin{eqnarray}
	\label{ppp}\lefteqn{
		\int_{\tilde{\phi}(\Omega)}     \frac{\partial \Delta^{h}v_1}{\partial y_i} \frac{\partial  \Delta^kv_2}{\partial y_j} \left(\frac{\partial\zeta _j}{\partial y_i}+\frac{\partial \zeta _i}{\partial y_j}\right)dy
				=- \int_{\partial \tilde{\phi}(\Omega)}  \Delta^{h}v_1 \Delta^{k+1}v_2 \zeta \cdot\nu d\sigma } \nonumber \\
& & \qquad		+\int_{\tilde{\phi}(\Omega)}\Delta^{h}v_1\nabla\Delta^{k+1}v_2\cdot\zeta  dy+\int_{\tilde{\phi}(\Omega)}\Delta^{h}v_1\Delta^{k+1}v_2{\rm div }\zeta  dy\nonumber \\
& & \qquad 		+\int_{\tilde{\phi}(\Omega)}\nabla \Delta^{h}v_1\cdot  \nabla\Delta^{k}v_2 {\rm div }\zeta  dy-\int_{\tilde{\phi}(\Omega)}\Delta^{h+1}v_1\nabla\Delta^kv_2\cdot\zeta  dy.
	\end{eqnarray}
	
By observing that the first summand in the right-hand side of (\ref{ppp}) vanish if $k< s-1$, and by combining (\ref{ppp0}) and (\ref{ppp}), we obtain a telescopic sum and we deduce the validity of (\ref{applem}). \hfill $\Box$

\begin{thm}
\label{beldif}
 Let $\Omega $ be a bounded open set in ${\mathbb{R}}^N$ of class $C^1$, $n\in {\mathbb{N}}$,   $u_1, u_2\in W^{n,2}_0(\Omega)$. Let $\tilde \phi \in {\mathcal{A}}_{\Omega}^{n}$ and  $v_i=u_i\circ \tilde \phi ^{(-1)}$, $i=1,2$. Assume that $\tilde \phi (\Omega )$ is of class $C^1$
and that  $v_1,v_2\in W^{2n,2}(\tilde \phi (\Omega))$. Then
\begin{eqnarray}\label{beldiff1}\lefteqn{
({\rm d|}_{\phi =\tilde\phi } (-\Delta _{\phi})^n[\psi])[u_1][u_2]
=-\int_{\partial \tilde \phi (\Omega )}\frac{\partial^nv_1}{\partial \nu^n}
\frac{\partial^nv_2}{\partial \nu^n}\zeta\cdot \nu d\sigma } \nonumber \\ & & \qquad\qquad\qquad\qquad\qquad\,  -\int _{ \tilde \phi (\Omega )} ((-\Delta )^{n}v_1\nabla v_2 +
(-\Delta )^{n}v_2\nabla v_1)\cdot  \zeta dy,
\end{eqnarray}
for all $\psi \in C^{n}_b(\Omega\, ;{\mathbb{R}}^N)$, where $\zeta =\psi \circ \tilde\phi^{(-1)}$ .
\end{thm}

{\bf Proof.} First, we consider the case where $n$ is an even number of the form $n=2s$ with $s\in {\mathbb{N}}_0$.  By formula (\ref{pol1}) and standard calculus we have
 \begin{eqnarray}
	\label{this}\lefteqn{
	\left(d|_{\phi=\tilde{\phi}}\Delta_{\phi}^{2s}[\psi](u_1)\right)(u_2)
	=\int_{\Omega}d|_{\phi=\tilde{\phi}}\Delta_{\phi}^{s}u_1[\psi]\Delta_{\tilde{\phi}}^su_2|\det D\tilde{\phi}|dx }\nonumber  \\ & & \qquad\qquad\qquad\qquad\qquad\qquad
	+\int_{\Omega}\Delta_{\tilde{\phi}}^su_1d|_{\phi=\tilde{\phi}}\Delta_{\phi}^{s}u_2[\psi]|\det D\tilde{\phi}|dx\nonumber \\ & & \qquad\qquad\qquad\qquad\qquad\qquad
	+\int_{\Omega}\Delta_{\tilde{\phi}}^su_1\Delta_{\tilde{\phi}}^su_2d|_{\phi=\tilde{\phi}}|\det D\phi|[\psi]dx	.
	\end{eqnarray}
	Moreover, by standard calculus
\begin{equation}
\label{der4}
\left[
\left({\rm d}|_{\phi =\tilde \phi } \left({\mathrm{det}}\nabla\phi\right)[\psi]\right)\circ
\tilde \phi^{(-1)}\right]{\mathrm{det}}\nabla\tilde \phi^{(-1)}=
{\mathrm{div}}\left(\psi\circ\tilde  \phi^{(-1)} \right)\, ,
\end{equation}
hence
\begin{equation}\label{der5}
\int_{\Omega}\Delta_{\tilde{\phi}}^su_1\Delta_{\tilde{\phi}}^su_2d|_{\phi=\tilde{\phi}}|\det D\phi|[\psi]dx=
\int_{\tilde\phi (\Omega )}\Delta^s \tilde v_1\Delta^s \tilde v_2 {\rm div }\zeta dy.
\end{equation}
Formula (\ref{beldiff1}) easily follows by combining formulas (\ref{applem}), (\ref{this}), (\ref{der5}), by integrating by parts and  by observing that
$\Delta ^{s}v_i=\frac{\partial ^{2s} v_i}{\partial \nu ^{2s}}  $ on $\partial\tilde \phi (\Omega )$ since $v_i\in W^{2s,2}_0(\tilde \phi (\Omega))$.

Now, we consider the case where $n$ is an odd number of the form  $n=2s+1$ with $s\in {\mathbb{N}}_0$.  By formula (\ref{pol2}) and standard calculus we have

\begin{eqnarray}
	\label{twoodd}
	\left(d|_{\phi=\tilde{\phi}}\Delta^{2s+1}_{\phi}[\psi]( u_1)\right)( u_2)
		=\int_{\tilde \phi (\Omega )}\nabla\Delta^s v_1(\nabla \zeta+\nabla^t \zeta)\nabla^t\Delta^s v_2dy\nonumber \\
	-\int_{\tilde \phi (\Omega )}  \nabla\Delta^s v_1\nabla\Delta^s v_2{\rm div }\zeta dy
			-\int_{\tilde \phi (\Omega )}\nabla\Delta^s v_1\nabla\left((d|_{\phi=\tilde{\phi}}\Delta^s_{\phi} u_2[\psi])\circ \tilde{\phi}^{(-1)}\right)dy\nonumber  \\
	-\int_{\tilde \phi (\Omega )} \nabla\left((d|_{\phi=\tilde{\phi}}\Delta^s_{\phi} u_1[\psi])\circ \tilde{\phi}^{(-1)}\right)\nabla\Delta^s v_2dy.
	\end{eqnarray}
Moreover, integrating by parts yields
	
\begin{eqnarray}\label{twoodd1}
\lefteqn{ \int_{\tilde \phi (\Omega )}\nabla\Delta^s v_1(\nabla\zeta+\nabla^t\zeta )\nabla^t\Delta^s v_2dy } \nonumber  \\ & &\qquad
 =\sum_{h,k=1}^N \int_{\tilde \phi (\Omega )}\left(\frac{\partial \zeta_h}{\partial  y_k}\frac{\partial  \Delta^s v_1}{\partial  y_h}\frac{\partial \Delta^s v_2}{\partial  y_k}+\frac{\partial  \zeta_k}{\partial  y_h}\frac{\partial \Delta^s v_1}{\partial y_h}\frac{\partial \Delta^s v_2}{\partial  y_k}\right)dx \nonumber  \\ & & \qquad
=2 \int_{\partial \tilde \phi (\Omega )}  \frac{\partial  \Delta^s v_1}{\partial  \nu}\frac{\partial  \Delta^s v_2}{\partial \nu}\zeta\cdot\nu d\sigma \nonumber  \\
& & \qquad -\sum_{h,k=1}^N\int_{\tilde \phi (\Omega )}\left(\frac{\partial^2  \Delta^s v_1}{\partial  y_h\partial  y_k}\frac{\partial \Delta^s v_2}{\partial y_k}\zeta_h+\frac{\partial \Delta^s v_1}{\partial y_h}\frac{\partial ^2\Delta^s v_2}{y_k^2}\zeta_h\right)dy\nonumber  \\
& & \qquad -\sum_{h,k=1}^N\int_{\tilde \phi (\Omega )}\left(\frac{\partial ^2\Delta^s v_1}{\partial y_h^2}\frac{\partial \Delta^s v_2}{\partial y_k}\zeta_k+\frac{\partial \Delta^s v_1}{\partial  y_h}\frac{\partial ^2\Delta^s v_2}{\partial  y_h\partial y_k}\zeta_k\right)dy\nonumber \\
& & \qquad =\int_{\partial \tilde \phi (\Omega )}  \frac{\partial \Delta^s v_1}{\partial \nu}\frac{\partial \Delta^s v_2}{\partial \nu}\zeta\cdot\nu d\sigma+\int_{\tilde \phi (\Omega )}\nabla\Delta^s v_1\nabla\Delta^s v_2{\rm div} \zeta dy\nonumber  \\
& & \qquad -\int_{\tilde \phi (\Omega )}(\Delta^{s+1} v_1\nabla\Delta^s v_2+\Delta^{s+1} v_2\nabla\Delta^s v_1)\cdot\zeta dy.
\end{eqnarray}
	
By integrating by parts, changing variables in integrals and using formula (\ref{applem}), we obtain
\begin{eqnarray}\label{twoodd2}\lefteqn{\int_{\tilde \phi (\Omega )}\nabla\Delta^s v_i\nabla\left((d|_{\phi=\tilde{\phi}}\Delta^s_{\phi} u_j[\psi])\circ \tilde{\phi}^{(-1)}\right)dy} \nonumber \\ & &\qquad\qquad = -\int_{\Omega }\Delta^{s+1}_{\tilde \phi } u_i d|_{\phi=\tilde{\phi}}\Delta^s_{\phi} u_j[\psi]|{\rm det}\nabla\tilde \phi | dx\nonumber \\ & &\qquad\qquad = -
\int_{\tilde \phi (\Omega )}(\Delta^{s+1}v_i\nabla\Delta^sv_j-\Delta ^{2s+1}v_i\nabla v_j)\cdot \zeta dy
\end{eqnarray}
for all $i,j\in \{1,2\}$. Finally, formula (\ref{beldiff1}) easily follows by combining formulas (\ref{twoodd})-(\ref{twoodd2}) and by observing that
$\frac{\partial \Delta ^{s}v_i}{\partial \nu }=\frac{\partial ^{2s+1} v_i}{\partial \nu ^{2s}}  $ on $\partial\tilde \phi (\Omega )$ since $v_i\in W^{2s+1,2}_0(\tilde \phi (\Omega))$.\hfill $\Box$\\

{\bf Acknowledgments}: The second author acknowledges financial support from the research project PRIN 2008 ``Aspetti geometrici delle equazioni alle derivate parziali e questioni connesse''.\\

\noindent {\small
Davide Buoso$^{1}$ and Pier Domenico Lamberti$^{2}$\\
Dipartimento di Matematica\\
Universit\`{a} degli Studi di Padova\\
Via Trieste,  63\\
35126 Padova\\
Italy\\
e-mail$^{1}$:	dbuoso@math.unipd.it \\
e-mail$^{2}$:	lamberti@math.unipd.it}

\end{document}